\documentclass{article}
\usepackage{amsmath}
\usepackage{amssymb}
\usepackage{graphicx}
\usepackage{epsfig}
\begin{document}

\title{On the Behavior of the Residuals in Conjugate Gradient Method} 
\author{Teruyoshi Washizawa \\SARD Center, Canon Inc.}
\maketitle

\begin{abstract}
In conjugate gradient method, it is well known that the recursively computed residual differs from true one as the iteration proceeds in finite arithmetic. Some work have been devoted to analyze this be-havior and to evaluate the lower and the upper bounds of the difference. This paper focuses on the behavior of these two kinds of residuals, especially their lower bounds caused by the loss of trailing digit, respectively.
\end{abstract}

\section{Introduction}
Conjugate gradient (CG) method and its varieties are popular as one of the best unsteady iterative methods for solving the following linear system:
\begin{equation}
A x = b
\label{linear_equation}
\end{equation}
In CG method, an approximate solution $x_{k}$ is expected to approach the exact solution x*. For the symmetric positive definite A, it is proved that the A-norm of the error monotonically decreases as the iteration proceeds in exact arithmetic.
This will be called as A-norm monotonicity of the error in the remaining part of this article.
It is obvious that we cannot calculate directly such a norm of the error without the solution.
Therefore, almost all algorithms employ the residual which is easily calculated as the difference between the left hand side (LHS) and the right hand side (RHS) of (1), $r_{k} := b - Ax_{k}$.
In practice, the residual is calculated by the recursion formula because of the computational complexity of the matrix vector product $Ax_{k}$ \cite{Ginsburg1963,Bollen1984}.
However, this recursion formula causes another problem in which the recursive residual differs from the true residual as the iteration proceeds.
It can be also observed that the recursive residual decreases after the true one seams to reach its lower bound.
We should terminate the CG steps just before the difference is too large to be neglected.
Ginsburg has proposed a simple criterion \cite{Ginsburg1963}:
\begin{quote}
For the true residual calculated as the difference between LHS and RHS of a linear system and the recursive residual calculated by using the recursion formula, the procedure is terminated when the 2-norm of their difference is greater than the 2-norm of the recursive residual:
\begin{equation}
\| r_{k} \| < exp(k/n)^{2} \| s_{k} - r_{k} \| \nonumber
\end{equation}
where $n$ is dimensionarity of a linear system.
\end{quote}

Several researchers have proposed the estimations of the lower and the upper bound of the norm of the error and the residual.
Wo\'{z}niakowski investigated the numerical stabilities and good-behaviors of three stationary iterative methods and CG method using the true residual $b-Ax_{k}$ \cite{Wozniakowski1978,Wozniakowski1980}.
Wo\'{z}niakowski gave the upper bound of the ultimately attainable accuracies of the $A$-norm and the 2-norm of the error, and 2-norm of the true residual.
Bollen gives the round-off error analysis of descent methods and lead a general result on the attainable accuracy of the approximate solution in finite arithmetic \cite{Bollen1984}.
It has also shown that the general result is applied to the Gauss-Southwell method and the gradient method to obtain the decreasing rates of the $A$-norm of the error in finite arithmetic.
Greenbaum have shown that for tiny perturbation $\epsilon_{M}$, the eigenvalues and the $A$-norm of the error vectors generated over a fixed number of perturbed itera-tive steps are approximately the same as those quantities generated by the exact recurrences applied to a "nearby" matrix \cite{Greenbaum1989}.
The lower bound of the true residual is pointed out in \cite{Greenbaum1997}.
Two kinds of the estimates of the A-norm of the error at every step in CG algorithm has been proposed and verified that those estimates are the lower and the upper bound in \cite{Meurant1997}.
The lower and the upper bounds of the A-norm of the error have been also given by Meurant\cite{Meurant1997} and Strako\v{s} and Tich\'{y}\cite{Strakos:Tichy2002}.
Strako\v{s} and Tich\'{y} have proposed the tight estimate for the lower bound of both the A-norm and the 2-norm of the error in every step.
This stepwise lower bound, however, keeps decreasing after the error reaches its 'global' lower bound.
Therefore, the terminating criteria by using this stepwise lower bound cannot detect the global lower bound of the error.
Calvetti et al. has proposed the estimates of the lower and upper bound of the A-norm of the error in CG method \cite{Calvetti:Morigi:Reichel:Sgallari2000}.
Those previous studies give the stepwise lower and the upper bound of the error and the residual but the global bounds.
In the remaining part of this article, we will first show that the true and the recursive residual almost monotonically de-crease as the iteration proceeds. Then, these lower bounds will be shown.

\section{Notation}
We shall give the notations appeared throughout this article.
$A$ and $b$ is, respectively, a coefficient matrix and a constant vector in a linear system.
$\| \cdot \|$ in connection with a vector and a matrix, respectively, stands for the 2-norm and spectral norm, $\| \cdot \|_{A}$ in connection with a vector stands for the norm under the metric tensor $A$.
The exact value of a variable $x$ is denoted as $\bar{x}$.
The floating point representation of a variable $x$ is denoted simply as $x$.
The computational error caused by the floating point representation is denoted as an operator $\epsilon_{M}(x) := x - \bar{x}$.
The exact solution of (1) is denoted as $x^{*}$ which is described formally as $x^{*} = A^{-1}b$.
At the k-th step, an approximate solution, the error, the true residual, and the recursive residual is, respectively, described as $x_{k}$, $e_{k}$, $s_{k}$, and $r_{k}$.
They are computed in CG method as follows:
\begin{eqnarray}
x_{k+1} & := & x_{k} + \alpha_{k}p_{k}, \nonumber \\
e_{k} & := & x^{*} - x_{k}, \nonumber \\
s_{k} & := & b - Ax_{k}, \nonumber \\
r_{k+1} & := & r_{k} - \alpha_{k}Ap_{k} \nonumber
\end{eqnarray}
Since $A$ is a constant matrix and $\|\epsilon_{M}(A)\|/\|A\|$ is almost equal to $\epsilon_{M}$ without the dependence on the number of iterations, $\epsilon_{M}(A)$ is out of our concern as well as $\epsilon_{M}(b)$.

\section{Almost monotonicity of residuals in finite arithmetic}
In this section, we will see the true and the recursive residual has the 2-norm almost monotonicity in finite arithmetic, respectively.

\subsection{The 2-norm almost monotonicity of true residual}
The true residual is calculated as the difference between LHS and RHS. This is equivalent to multiplication of A to the error in finite arithmetic,
\begin{equation}
s_{k} = b - Ax_{k} = b - A(x_{k} - x^{*}) - Ax^{*} = - Ae_{k} \nonumber
\end{equation}
The behavior of true residual $s_{k}$ is, therefore, equivalent to $Ae_{k}$.
The $A$-norm monotonicity of the error in finite arithmetic has been proved in theorem-3.1 of [2].
The following theorem shows the error has the 2-norm almost monotonicity.

\begin{quote}
{\bf Theorem-1.}
If $e_{k}$ has the $A$-norm monotonicity, $e_{k}$ has the 2-norm almost monotonicity for a regular matrix $A$,
\begin{equation}
\exists k > j, \| e_{k} \| < \| e_{j} \|
\label{theorem01}
\end{equation}
\\{\bf Proof.}
The relationship between 2-norm and $A$-norm of an error is
\begin{equation}
\| e_{k} \| = \| A^{-1/2}e_{k} \|_{A} \le \|A^{-1/2} \| \| e_{k} \|_{A}
\label{proof01_01}
\end{equation}
similarly,
\begin{equation}
\| e_{k} \|_{A} = \| A^{1/2}e_{k} \| \le \|A^{1/2} \| \| e_{k} \|
\label{proof01_02}
\end{equation}
Substituting (\ref{proof01_01}) and (\ref{proof01_02}) into (\ref{theorem01}) and we yield
\begin{equation}
\| A^{-1/2} \| \| e_{k} \|_{A} < \|A^{1/2} \|^{-1} \| e_{k} \|_{A}
\label{proof01_03}
\end{equation}
Equation (\ref{proof01_03}) holds if $k$ exists to satisfy the following relation
\begin{equation}
\| e_{k} \|_{A} < \kappa(A^{1/2})^{-1} \| e_{k} \|_{A}
\label{proof01_04}
\end{equation}
where $\kappa(A)$ is the condition number of a matrix $A$.
From the $A$-norm monotonicity of the error $e_{k}$, since the following relation holds for any positive value $a$,
\begin{equation}
\exists k > j, \| e_{k} \|_{A} < a \| e_{j} \|_{A} \nonumber
\end{equation}
there exists $k>j$ satisties (\ref{proof01_04}) and consequently (\ref{theorem01}).
We have to notice that (\ref{theorem01}) does not hold when $a=0$, i.e., the inverse of the coefficient matrix $A$ is singular.
\end{quote}
Theorem-1 leads to the 2-norm almost monotonicity of the true residual using the relationship $\| s_{k} \| = \| Ae_{k} \|$.

\begin{quote}
{\bf Theorem-2.}
If $e_{k}$ has the 2-norm almost monotonicity, $s_{k}$ has the 2-norm almost monotonicity for a regular matrix $A$,
\begin{equation}
\exists k > j, \| s_{k} \| < \| s_{j} \|
\label{theorem02}
\end{equation}
\\{\bf Proof.}
The relationship between the 2-norm of the error and that of the residual is
\begin{equation}
\| s_{k} \| = \| Ae_{k} \| \le \| A \| \| e_{k} \|
\label{proof02_01}
\end{equation}
Similarly,
\begin{equation}
\| e_{j} \| = \| A^{-1}s_{j} \| \le \| A^{-1} \| \| s_{j} \|
\label{proof02_02}
\end{equation}
Substituting (\ref{proof02_01}) and (\ref{proof02_02}) into (\ref{theorem02}) and we yield
\begin{equation}
\| A \| \| e_{k} \| < \| A^{-1} \|^{-1} \| e_{j} \|
\label{proof02_03}
\end{equation}
Equation (\ref{proof02_03}) holds if $k$ exists to satisfy the following relationship
\begin{equation}
\| e_{k} \| < \kappa(A)^{-1} \| e_{j} \|
\label{proof02_04}
\end{equation}
From the 2-norm almost monotonicity of the error $e_{k}$, since the following relationship holds for any positive value a,
\begin{equation}
\exists k > j, \| e_{k} \| < a \| e_{j} \|
\label{proof02_05}
\end{equation}
there exists $k>j$ that satisfies (\ref{proof02_04}) and consequently (\ref{theorem02}).
\end{quote}

\section{The 2-norm almost monotonicity of recursive residual}
Before the proof of almost monotonicity of recursive residual in finite arithmetic, we first give the proof of almost monotonicity of recursive residual in exact arithmetic.
From the $A$-norm monotonicity of the error in exact arithmetic, the 2-norm almost monotonicity of the residual in exact arithmetic can be proved.
We have to notice that the recursive residual $\bar{r}_{j}$ is identical to the true residual $\bar{s}_{j}$ in exact arithmetic.

\begin{quote}
{\bf Theorem-3.}
If $\forall n, \| \bar{e}_{n+1} \|_{A} < \| \bar{e}_{n} \|_{A}$, then the following propositoin holds for a regular matrix A:
\begin{equation}
\exists k > j, \| \bar{r}_{k} \| < \| \bar{r}_{j} \|
\label{theorem03}
\end{equation}
\\{\bf Proof.}
The relationship between the error and the residual gives:
\begin{equation}
\| \bar{r}_{k} \| = \| A \bar{e}_{k} \| = \| A^{1/2} \bar{e}_{k} \|_{A}
\label{proof03_01}
\end{equation}
Then we yield the lower and the upper bound of the 2-norm of the true residual:
\begin{equation}
\| A^{-1/2} \|^{-1} \| \bar{e}_{k} \|_{A} \le \| \bar{r}_{k} \| \le \| A^{1/2} \| \| \bar{e}_{k} \|_{A}
\label{proof03_02}
\end{equation}
From above equation, the sufficient condition for (\ref{theorem03}) can be given as follows:
\begin{equation}
\| A^{1/2} \| \| \bar{e}_{k} \|_{A} < \| A^{-1/2} \|^{-1} \| \bar{e}_{j} \|_{A}
\label{proof03_03}
\end{equation}
that is, the equation holds if there exists $k>j$ so that
\begin{equation}
\| \bar{e}_{k} \|_{A} < \kappa(A^{1/2})^{-1} \| \bar{e}_{j} \|_{A}
\label{proof03_04}
\end{equation}
From the $A$-norm monotonicity of the error $\bar{e}_{k}$, the following equation holds for any positive value $a$,
\begin{equation}
\exists k > j, \| \bar{e}_{k} \|_{A} < a \| \bar{e}_{j} \|_{A}
\label{proof03_05}
\end{equation}
and (\ref{theorem03}) holds.
\end{quote}
Now we show the almost monotonicity of the recursive residual in finite arithmetic.
\begin{quote}
{\bf Theorem-4.}
If the recursive residual has 2-norm almost monotonicity in exact arithmetic, then the recursive residual has the 2-norm almost monotonicity in finite arithmetic:
\begin{equation}
\exists k > j, \| r_{k} \| < \| r_{j} \|
\label{theorem04}
\end{equation}
\\{\bf Proof.}
Equation (19) is rewritten as
\begin{equation}
\exists k > j, \| \bar{r}_{k} + \epsilon_{M}(\bar{r}_{k}) \| < \| \bar{r}_{j} + \epsilon_{M}(\bar{r}_{j}) \| \nonumber
\end{equation}
The following relationship is one of its sufficient conditions
\begin{equation}
\max \left[ \| \bar{r}_{k} + \epsilon_{M}(\bar{r}_{k}) \| \right] < \min \left[ \| \bar{r}_{j} + \epsilon_{M}(\bar{r}_{j}) \| \right]
\label{proof04_01}
\end{equation}
The evaluation of the maximum value of LHS is
\begin{equation}
\max \left[ \| r_{k} \| \right] = (1 + \epsilon_{M}) \| \bar{r}_{k} \|
\label{proof04_02}
\end{equation}
Similarly, the minimum value of RHS is evaluated as
\begin{equation}
\min \left[ \| r_{j} \| \right] = (1 - \epsilon_{M}) \| \bar{r}_{j} \|
\label{proof04_03}
\end{equation}
Substituting (\ref{proof04_02}) and (\ref{proof04_03}) into (\ref{proof04_01}), the sufficient condition (20) is given as
\begin{equation}
\exists k > j, \| \bar{r}_{k} \| < (1 - \epsilon_{M}) / (1 + \epsilon_{M}) \| \bar{r}_{j} \|
\label{proof04_04}
\end{equation}
There exists $k>j$ for $(1 - \epsilon_{M}) / (1 + \epsilon_{M}) > 0$ from theorem-3 and (\ref{proof04_04}) holds.
\end{quote}

\section{Lower bounds of error and residual in finite arithmetic}
It has been shown that the 2-norm of two kinds of residuals, respectively, decreases almost monotonically in finite arithmetic in the previous section.
Now we consider whether if the 2-norm of each variable stops decreasing before the approximates $x_{k}$ does not reach its target $x^{*}$.
\subsection{Lower bound of error}
Theoem-1 shows the approximate solution $x_{j}$ approaches the exact solution almost monotonically in finite arithmetic.
The correction of the recursion formula of $x_{k}$, however, vanishes by the loss of trailing digits so that the error stops changing, i.e.,
\begin{equation}
| \Delta x_{stop}(n)| / |x_{stop}(n)| < \epsilon_{M} \Longrightarrow x_{stop+1} = x_{stop} \nonumber
\end{equation}
where $x_{k}(n)$ is the n-th component of $x_{k}$.

On the other hand, the solution in finite arithmetic $x^{*}$ is not always identical to that in exact arithmetic  $\bar{x}^{*}$.
Therefore, the target for iterative algorithms in finite arithmetic should not be $\bar{x}^{*}$ but $x^{*}$.
The error caused by the loss of trailing digits is described formally as $x^{*} - x_{stop}$.
Since the true residual is given by multiplying A to the error, the lower bound of the true residual is given as $A(x^{*} - x_{stop})$.

\subsection{Lower bound of recursive residual}
Theorem-4 shows the recursive residual reduces its 2-norm almost monotonically. The next theorem proves that the change of the recursive residual stops only when $\|r_{k+1} \| < \epsilon_{M} \| r_{k} \|$.
It decreases almost monotonically until then.
\begin{quote}
{\bf Theorem-5.}
The recursive residual never have a lower bound caused by the loss of trailing digits.
\\{\bf Proof.}
The recursion formula of the residual is in general described as
\begin{equation}
r_{k+1} = r_{k} - \Delta k_{k}
\label{proof05_01}
\end{equation}
where $\Delta kr_{k} := \alpha_{k}Ap_{k}$.
The residual reaches its lower bound $r_{k}$ if the following condition is satisfied:
\begin{equation}
\Delta r_{k} < \epsilon_{M}(r_{k})
\label{proof05_02}
\end{equation}
We will show the condition of (\ref{proof05_02}) never be satisfied in not only exact but also finite arithmetic.
In exact arithmetic, (\ref{proof05_01}) satisfies the following relationship:
\begin{eqnarray}
\| r_{k+1} \|^{2} & = & \| \bar{r}_{k} - \Delta \bar{r}_{k} \|^{2} = \| \bar{r}_{k} \|^{2} - 2(\bar{r}_{k}, \Delta \bar{r}_{k}) + \| \Delta \bar{r}_{k} \|^{2} \nonumber \\
& = & \| \bar{r}_{k} \|^{2} - 2(\bar{r}_{k}, \bar{r}_{k} - \bar{r}_{k+1}) + \| \Delta \bar{r}_{k} \|^{2} \nonumber \\
& = & - \| \bar{r}_{k} \|^{2} + \| \Delta \bar{r}_{k} \|^{2} \nonumber
\end{eqnarray}
where using $(\bar{r}_{k}, \bar{r}_{k+1})=0$.
The following relationship holds directly from above equation in exact arithmetic :
\begin{equation}
\| \Delta \bar{r}_{k} \| \ge \| \bar{r}_{k} \|
\label{proof05_03}
\end{equation}
This shows that (\ref{proof05_02}) never be satisfied in exact arithmetic.
Now we evaluate above in finite arithmetic.
\begin{eqnarray}
\frac{\| \Delta r_{k} \|}{\| r_{k} \|} & \ge & \frac{\min \left[ \Delta \bar{r}_{k} + \epsilon_{M}(\Delta \bar{r}_{k}) \right]}{\max \left[ \| \bar{r}_{k} + \epsilon_{M}(\bar{r}_{k} ) \| \right]} \nonumber \\
& = & \frac{(1 - \epsilon_{M})\| \Delta \bar{r}_{k} \|}{(1 + \epsilon_{M})\| \bar{r}_{k} \|} = \left( 1 - \frac{2 \epsilon_{M}}{1 + \epsilon_{M}} \right) \frac{\| \Delta \bar{r}_{k} \|}{\| \bar{r}_{k} \|} \nonumber
\end{eqnarray}
According to (\ref{proof05_03}), we yield
\begin{equation}
\frac{\| \Delta r_{k} \|}{\| r_{k} \|} \ge 1 - \frac{2 \epsilon_{M}}{1 + \epsilon_{M}} \nonumber
\end{equation}
and prove that the recursive residual never have a lower bound caused by the loss of trailing digits.
The termination of the iterations is caused only when $\| r_{k+1} \| \le \epsilon_{M} \| r_{k} \|$ by the significant decrease of the recursive residual for $\alpha_{k} Ap_{k} \approx r_{k}$.
\end{quote}

\section{Conclusions}
In this article, the convergence behaviors of true and recursive residual have been analyzed.
The results obtained are summarized below:
\begin{itemize}
\item In finite arithmetic, the 2-norm of the error and the residual, respectively, almost monotonically decreases.
\item 2-norm of the error has the lower bound in finite arithmetic as well as the true residual.
\item 2-norm of the recursive residual never have a non-zero lower bound caused by the loss of trailing digits in finite arithmetic.
\end{itemize}


\begin{thebibliography}{99}
\bibitem{Ginsburg1963}
T. Ginsburg : "The conjugate gradient method," Numerische Mathematik, 5(1), pp.191-200 (1963).
\bibitem{Bollen1984}
J.A.M. Bollen : gNumerical Stability of Descenet Methods for Solving Linear Equations,h Numerische Mathematik, 43, pp.361-377, 1984.
\bibitem{Wozniakowski1978}
H. Wo\'{z}niakowski :gRound-Off Error Analysis of Iteratinos for Large Linear Systems,h Numeriche Mathematik, 30, pp.301-314, 1978.
\bibitem{Wozniakowski1980}
H. Wo\'{z}niakowski :gRoundoff Error Analysis of New Class of Conjugate-Gradient Algorithms,h Linear Algebra and its Applications, 29, pp.507-529, 1980.
\bibitem{Greenbaum1989}
A. Greenbaum : "Behavior of slightly perturbed Lanczos and conjugate-gradient recurrences," Linear Algebra and Its Applications, 113, pp.7-63 (1989).
\bibitem{Greenbaum1997}
A. Greenbaum : "Estimating the Attainable Accuracy of Recursively Computed Residual Methods,h SIAM Journal on Matrix Analysis and Applications, 18(3), pp.535-551, 1997.
\bibitem{Meurant1997}
G. Meurant, "The computation of bounds for the norm of the error in the conjugate gradient algorithm," Numerical Algorithms, 16, pp.77-87, 1997.
\bibitem{Strakos:Tichy2002}
Z. Strako\v{s}, and P. Tich\'{y}, "On error estimation in the conjugate gradient method and why it works in finite precision computations," Electronic Transactions on Numerical Analysis, Vol.13, pp.56-80, 2002.
\bibitem{Calvetti:Morigi:Reichel:Sgallari2000}
D. Calvetti, S. Morigi, L. Reichel, and F. Sgallari, "Computable error bounds and estimates for the conjugate gradient method," Numerical Algorithms, 25, pp.75-88, 2000.
\end{thebibliography}
\end{document}